\documentclass[12pt]{article}
\usepackage{amssymb}

\newtheorem{thm}     {Theorem}[section]

\newtheorem{definition}  [thm]{Definition}

\newtheorem{lemma}   [thm]{Lemma}

\newcommand{\proof} {\noindent{\bf Proof. }}

\newcommand{\B}{\mathbb B}
\newcommand{\C}{\mathbb C}
\newcommand{\D}{\mathbb D}

\newcommand{\R}{\mathbb R}

\newcommand{\st}{{\rm st}}

\def\bar{\overline}

\begin{document}

\title{On boundary properties of asymptotically holomorphic functions}
\author{Alexandre Sukhov{*} }
\date{}
\maketitle

{\small

* Universit\'e  de Lille, Laboratoire
Paul Painlev\'e,
Departement de
Math\'ematique, 59655 Villeneuve d'Ascq, Cedex, France, sukhov@math.univ-lille1.fr
The author is partially suported by Labex CEMPI.

Institut of Mathematics with Computing Centre - Subdivision of the Ufa Research Centre of Russian
Academy of Sciences, 45008, Chernyshevsky Str. 112, Ufa, Russia.

}
\bigskip

{\small Abstract.  We prove a   Fatou type theorem for   bounded functions with  
$\overline\partial_J$ differential of a controled growth on smoothly bounded domains in an almost complex manifold.}



\bigskip


\section{Introduction}

This paper is a continuation of the work \cite{Su2}. We improve the main results of \cite{Su2} establishing a general version of the Chirka - Lindel\"of principle and the Fatou type  theorem for bounded asymptotically holomorphic functions on almost complex manifolds. These functions admit the antiholomorphic part of the differential satisfying some asymptotic growth conditions near the boundary. Such classes of functions naturally appear in Several Complex Variables, PDE and related topics. Our results extend the known results 
\cite{Ch,Fo,Sa1,St} obtained for the case of $\C^n$ with the standard complex structure. Note that  in this case our results also are knew.

\section{ Almost complex manifolds  and almost holomorphic functions} This is a preliminary section. We recall basic notions of the almost complex geometry making the  presentation of our results more convenient. 
Everywhere through this paper we assume that manifolds and almost complex structures are of class $C^\infty$ (the word "smooth" means the regularity of this class). However,  our main results are also valid under considerably weaker regularity assumptions.

\subsection{ Almost complex manifolds} Let $M$ be a smooth  manifold of real dimension $2n$. {\it An almost complex structure} $J$ on $M$ is a smooth map  which associates to every point $p \in M$ a linear isomorphism $J(p): T_pM \to T_pM$ of the tangent space $T_pM$ such that $J(p)^2 = -I_p$; here  $I_p$ denotes the identity map of $T_pM$. Thus, every linear operator $J(p)$ is a  complex structure on a vector space $T_pM$ in the usual sense of Linear Algebra. When $J$ is fixed, a couple $(M,J)$ is called {\it an almost complex manifold} of complex dimension n. 

A  fundamental example of an almost complex structure  is given by the {\it standard complex structure} $J_{st} = J_{st}^{(2)}$ on $M = \R^2$.  This linear operator is represented in the canonical coordinates of $\R^2$ by the matrix 

\begin{eqnarray}
\label{J_st}
J_{st}= \left(
\begin{array}{cll}
0 & & -1\\
1 & & 0
\end{array}
\right)
\end{eqnarray}
More generally, {\it the standard complex structure $J_{st}$ on $\R^{2n}$} is represented by the block diagonal matrix $diag(J_{st}^{(2)},...,J_{st}^{(2)})$ (usually we drop the notation of dimension because its value  will be clear from the context). Setting $iv := J_{st}v$ for a vector $v \in \R^{2n}$, we identify the real space $(\R^{2n},J_{st})$ with the complex linear space  $\C^n$; we  use the notation  $z = x + iy = x + J_{st}y$ for the standard complex coordinates $z = (z_1,...,z_n) \in \C^n$.

Let $(M,J)$ and $(M',J')$ be smooth  almost complex manifolds. A $C^1$-map $f:M' \to M$ is called  
{\it $(J’,J)$-complex or  $(J',J)$-holomorphic}  if it satisfies {\it the Cauchy-Riemann equations} 
\begin{eqnarray}
\label{CRglobal}
df \circ J' = J \circ df.
\end{eqnarray}

In particular   a map $f: \C^n \to \C^m$ is $(J_{st},J_{st})$-holomorphic if and only if each component of $f$ is a usual holomorphic function.

 Every almost complex manifold
$(M,J)$ can be viewed locally as the unit ball $\B$ in
$\C^n$ equipped with a small (in any $C^m$-norm) almost complex
deformation of $J_{st}$. The following well-known statement is often  useful.
\begin{lemma}
\label{lemma1}
Let $(M,J)$ be an almost complex manifold. Then for every point $p \in
M$, every  $m \geq 0$ and   $\lambda_0 > 0$ there exist a neighborhood $U$ of $p$ and a
coordinate diffeomorphism $z: U \rightarrow \B$ such that
$z(p) = 0$, $dz(p) \circ J(p) \circ dz^{-1}(0) = J_{st}$,  and the
direct image $ z_*(J): = dz \circ J \circ dz^{-1}$ satisfies $\vert\vert z_*(J) - J_{st}
\vert\vert_{C^m(\bar {\B})} \leq \lambda_0$.
\end{lemma}
\proof There exists a diffeomorphism $z$ from a neighborhood $U'$ of
$p \in M$ onto $\B$ satisfying $z(p) = 0$; after an additional linear change of coordinates one can achieve  $dz(p) \circ J(p)
\circ dz^{-1}(0) = J_{st}$ (this is a classical fact from the Linear Algebra). For $\lambda > 0$ consider the isotropic dilation
$h_{\lambda}: t \mapsto \lambda^{-1}t$ in $\R^{2n}$ and the composition
$z_{\lambda} = h_{\lambda} \circ z$. Then $\lim_{\lambda \rightarrow
0} \vert\vert (z_{\lambda})_{*}(J) - J_{st} \vert\vert_{C^m(\bar
{\B})} = 0$ for every  $m \geq 0$. Setting $U = z^{-1}_{\lambda}(\B)$ for
$\lambda > 0$ small enough, we obtain the desired statement. 
In what follows we often denote the structure $z_*(J)$ again by $J$ viewing it as a local representation of $J$ in the coordinate system $(z)$.

Recall that an almost complex structure $J$ is called {\it integrable} if $(M,J)$ is locally biholomorphic in a neighborhood of each point to an open subset of $(\C^n,J_{st})$. In the case of complex dimension 1 every almost complex structure is integrable. In the case of complex dimension $> 1$  integrable almost complex structures form a highly special subclass in the space of all almost complex structures on $M$; an efficient criterion of integrablity is provided by the classical theorem of Newlander - Nirenberg \cite{NeNi}.

\bigskip

 \subsection{Pseudoholomorphic discs} Let $(M,J)$ be an almost complex manifold of dimension $n > 1$. For a "generic" choice of an almost complex structure, any holomorphic (even locally) function on $M$ is constant because the Cauchy-Riemann equations are overdetermined. For the same reason $M$ does not admit non-trivial $J$-complex submanifolds  of complex dimension $> 1$. The unique  exceptional   case arises when  $J$-complex submanifolds are of  complex dimension 1. They always exist at least locally.

{\it Pseudoholomorphic curves} are parametrized  by the   solutions  $f$ of (\ref{CRglobal}) in the special case   where $M'$ has the complex dimension 1. These holomorphic maps are called $J$-complex (or $J$-holomorphic or {\it pseudoholomorphic} ) curves. Note that we view here the curves as maps i.e. we consider parametrized curves.
We use the notation  $\D = \{ \zeta \in \C: \vert \zeta \vert < 1 \}$ for  the
unit disc in $\C$ always assuming that it is equipped with the standard complex structure   $J_{\st}$. If in the equations (\ref{CRglobal})  we have $M' = \D$, we  call such a map $f$ a $J$-{\it complex  disc} or a  {\it pseudoholomorphic disc} or just a  holomorphic disc
when the structure  $J$ is fixed. 

A fundamental fact is that  pseudoholomorphic discs always exist in a suitable neighborhood of any point of $M$; this is the classical Nijenhuis-Woolf theorem (see \cite{NiWo}). Here it is convenient to rewrite the equations (\ref{CRglobal}) in local coordinates  similarly to the complex version of the usual Cauchy-Riemann equations.

 Everything will be local, so (as above) we are in a neighborhood $\Omega$ of $0$ in $\C^n$ with the standard complex coordinates $z = (z_1,...,z_n)$. We assume that $J$ is an almost complex structure defined on $\Omega$ and $J(0) = J_{st}$. Let 
$$z:\D \to \Omega,$$ 
$$z : \zeta \mapsto z(\zeta)$$ 
be a $J$-complex disc. Setting $\zeta = \xi + i\eta$ we write (\ref{CRglobal}) in the form $z_\eta = J(z) z_\xi$. This equation can be written as

\begin{eqnarray}
\label{holomorphy}
z_{\bar\zeta} - A(z)\bar z_{\bar\zeta} = 0,\quad
\zeta\in\D.
\end{eqnarray}
Here a smooth map $A: \Omega \to Mat(n,\C)$ is defined by the equality $L(z) v = A \overline v$ for any vector $v \in \C^n$ and $L$ is an $\R$-linear map defined by $L = (J_{st} + J)^{-1}(J_{st} - J)$. It is easy to check that the condition $J^2 = -Id$ is equivalent to the fact that $L$ is $\overline\C$-linear. The matrix $A(z)$ is called {\it the complex matrix} of $J$ in the local coordinates $z$. Locally the correspondence between $A$ and $J$ is one-to-one. Note that the condition $J(0) = J_{st}$ means that $A(0) = 0$. 

If $t$ are other local coordinates and $A'$ is the corresponding complex matrix of $J$ in the coordinates $t$, then, as it is easy to check, we have the following transformation rule:

\begin{eqnarray}
\label{CompMat}
A’ = (t_z A  + { t}_{\overline z})({\overline t}_{\overline z} + {\overline t}_{ z}A)^{-1}
\end{eqnarray}
(see \cite{SuTu}).

Note that one can view the equations (\ref{holomorphy}) as a quasilinear analog of the Beltrami equation for vector-functions. From this point of view, the theory of pseudoholomorphic curves is an analog of the theory of quasi-conformal mappings.

\bigskip

Recall that for a complex function $f$  {\it the Cauchy-Green transform} is defined by

\begin{eqnarray}
\label{CauchyGreen}
Tf(\zeta) = \frac{1}{2 \pi i} \int_{\D} \frac{f(\omega)d\omega \wedge d\overline\omega}{\omega - \zeta}
\end{eqnarray}
 This  classical linear integral operator  has the following properties (see \cite{Ve}):
\begin{itemize}
\item[(i)] $T: C^r(\D) \to C^{r+1}(\D)$ is a bounded linear operator for every non-integer $r > 0$ ( a similar property holds in the Sobolev scales, see below). Here we use the usual H\"older norm on the space $C^r(\D)$.
\item[(ii)] $(Tf)_{\overline\zeta} = f$ i.e. $T$ solves the $\overline\partial$-equation in the unit disc. 
\item[(iii)] the function $Tf$ is holomorphic on $\C \setminus \overline\D$.
\end{itemize}
Fix a real non-integer $r > 1$. Let $z: \D \to \C^n$, $z: \D \ni \zeta \mapsto z(\zeta)$ be a $J$-complex disc. 
Since  the operator
$$\Psi_{J}: z \longrightarrow w =  z - TA(z) \overline {z}_{\overline \zeta} $$
takes the space   $C^{r}(\overline{\mathbb D})$  into itself,  we can write   the
equation (\ref{CRglobal}) in the form 
$(\Psi_J(z))_{\overline \zeta}  = 0$. Thus, the disc $z$ is $J$-holomorphic if
and only if the map $\Psi_{J}(z):\mathbb D \longrightarrow \C^n$ is
$J_{st}$-holomorphic.
When the norm of $A$  is small enough (which is assured  by Lemma \ref{lemma1}),
then  by the implicit function theorem the operator    $\Psi_J$
is invertible  in the space $C^r(\D)$ and we obtain a bijective
correspondence between  $J$-holomorphic discs and usual
holomorphic discs. This easily implies the existence of a $J$-holomorphic disc
in a given tangent direction through a given point of $M$, as well as  a smooth dependence of such a
disc  on a deformation of a point or a tangent vector, or on an almost complex structure; this also establishes  the interior elliptic regularity of discs. This is the classical Nijenhuis-Woolf theorem, see \cite{NiWo}.

\bigskip

Let $(M,J)$ be an almost complex manifold and $E \subset M$ be a real submanifold of $M$. 
Suppose that a $J$-complex disc $f:\D \to M$ is  continuous on $\overline\D$.  With some abuse of terminology, we also call the image $f(\D)$  simply by a disc and we call the image $f(b\D)$ by  the boundary of a disc. If  $f(b\D) \subset E$, then we say that (the boundary of ) the disc  $f$ is {\it glued} or {\it attached} to $E$ or simply 
that $f$ is attached to $E$. If $\gamma \subset b\D$ is an arc and $f(\gamma) \subset E$, we say that $f$ is {\it glued or attached } to $E$ along $\gamma$.

\subsection{The $\overline\partial_J$-operator on an almost complex manifold $(M,J)$}

Consider now the second special class (together with pseudoholomorphic curves) of holomorphic maps. Consider first the situation when $J$ be an almost complex structure defined 
in a domain $\Omega\subset\C^n$; one can view this as a local coordinate representation of $J$ in a chart on $M$.

A $C^1$ function $F:\Omega\to\C$ is $(J,J_{st})$-holomorphic
if and only if it satisfies the Cauchy-Riemann equations
\begin{eqnarray}
\label{CRscalar}
F_{\bar z} + F_z A(z)  =0,
\end{eqnarray}
where $F_{\bar z} = (\partial F/\partial \overline{z}_1,...,\partial F/\partial \overline{z}_n)$ and $F_z = (\partial F/\partial {z}_1,...,\partial F/\partial {z}_n)$ are viewed as  row-vectors. Indeed, $F$ is $(J,J_{st})$ holomorphic if and only if 
for every $J$-holomorphic disc $z:\D \to \Omega$ the composition $F \circ z$ is a usual holomorphic function  that is $\partial (F \circ z) /\partial\overline{\zeta} = 0$ on $\D$. Then the  Chain rule in combination with (\ref{holomorphy}) leads to (\ref{CRscalar}). Generally the only solutions to (\ref{CRscalar}) are constant functions  unless $J$ is integrable (then $A$ vanishes identically in suitable coordinates). Note also that (\ref{CRscalar}) is a linear PDE system while (\ref{holomorphy}) is a quasilinear PDE system for a vector function on $\D$.

Every $1$-differential form $\phi$ on $(M,J)$ admits a unique decomposition $\phi = \phi^{1,0} + \phi^{0,1}$ with respect to $J$. In particular, if $F:(M,J) \to \C$ is a $C^1$-complex function, we have $dF = dF^{1,0} + dF^{0,1}$. We use the notation 
\begin{eqnarray}
\label{d-bar}
\partial_J F = dF^{1,0} \,\,\,\mbox{and}\,\,\, \overline\partial_J F = dF^{0,1}
\end{eqnarray}

In order to write these operators explicitely in local coordinates, we find a  local basis in the space of (1,0) and (0,1) forms. We view  $dz = (dz_1,...,dz_n)^t$ and $d\overline{z} = (d\overline{z}_1,...,d\overline{z}_n)^t$ as vector-columns. Then the forms 
\begin{eqnarray}
\label{FormBasis}
\alpha = (\alpha_1,..., \alpha_n)^t = dz - A d\overline{z} \,\,\, \mbox{and} \,\, \overline\alpha = d\overline{z} - \overline A dz
\end{eqnarray}
form a basis in the space of  (1,0) and (0,1) forms respectively. Indeed, it suffices to observe that a 1-form $\beta$ is of type (1,0) (resp. $(0,1)$)  if and only if for every $J$-holomorphic disc $z$ the pull-back $z^*\beta$ is a usual (1,0) (resp. $(0,1)$) form on $\D$. Using the equations (\ref{holomorphy}) we obtain the claim.

Now we decompose the differential $dF = F_zdz + F_{\overline{z}} d\overline{z} = \partial_J F + \overline\partial_J F$ with respect to  the basis $\alpha$, $\overline\alpha$ using (\ref{FormBasis}). We obtain the explicit expression 
\begin{eqnarray}
\label{d-bar2}
\overline\partial_J F = (F_{\overline{z}} (I - \overline{A}A)^{-1} + F_z (I - A\overline{A})^{-1}A)\overline\alpha
\end{eqnarray}

It is easy to check that the holomorphy condition $\overline\partial_J F = 0$ is equivalent to (\ref{CRscalar}) because $(I - A\overline{A})^{-1} A (I - \overline{A} A) = A$. Thus 

\begin{eqnarray*}
\overline\partial_J F = (F_{\overline{z}}  + F_z A)(I - \overline{A}A)^{-1}\overline\alpha
\end{eqnarray*}

We note that the matrix factor  $(I - A\overline A)^{-1}$ as well as the forms $\alpha$ affect only the non-essential constants in local estimates of the $\overline\partial_J$-operator near a boundary point which we will perfom in the next sections. So the reader  can assume that this operator is simply given by the left hand expression of (\ref{CRscalar}).

\begin{definition}
\label{DefSub}
Let $F$ be a complex function of class $C^1$ on a (bounded) domain $\Omega$ in an almost complex manifold $(M,J)$ of dimension $n$. We call a function $F$ a {\it subsolution of the $\overline\partial_J$ operator }  or simply 
$\overline\partial_J$-subsolution  on $\Omega$ if   is there exists  constant $C > 0$ and $\tau > 0$ such that
\begin{eqnarray}
\label{AlHol}
\parallel \overline\partial_J F (z)\parallel \le C dist(z,b\Omega)^{-1/2 + \tau}
\end{eqnarray}
for all $z \in \Omega$. Here we use the norm with respect to any fixed Riemannian metric on $M$.
\end{definition}

Obviously, non-constant $\overline\partial_J$-subsolutions exist in a sufficiently small neighborhhod of any point of $M$. For example  any function $F$ of class $C^1$ in an open neighborhhod of the compact set $\overline\Omega$ is a $\overline\partial_J$-subsolution on $\Omega$. Of course, any $C^1$ function $F$ with uniformly bounded  $\overline\partial_J F$ on $\Omega$,  satisfies (\ref{AlHol}). This subclass of functions was studied in \cite{Su2}. In the case of $\C^n$, a similar class of functions appeared in \cite{Fo}.

Let $F$ be a $\overline\partial_J$-subsolution on $\Omega$. Suppose that $A$ is the complex matrix of $J$ in a local chart $U$ and $z:\D \to U$ is a $J$-complex disc. It follows by the Chain Rule and (\ref{holomorphy}) that
$$(F \circ z)_{\overline\zeta} = (F_{\overline z} + F_zA){\overline z}_{\overline\zeta}.$$
Thus, if $h: \D \to \Omega$ is a $J$-complex disc of class $C^1(\overline\D)$, then the composition 
$F \circ h$ has a  $\overline\partial$-derivative satisfying (\ref{AlHol}) on $\D$ that is $F \circ h$ is a $\overline\partial_{J_{st}}$-subsolution on $\D$. Note that the constant $C$ and $\tau$ appearing in the upper bound of type (\ref{AlHol}) for the $\overline\partial (F \circ h)$, depend only on constants from the upper bound on $\overline\partial_J F$ in (\ref{AlHol}),  and the $C^1$ norm of $h$ on $\overline\D$ as well. In particular, if $(h_t)$ is a family of $J$-complex discs in $\Omega$ and $C^1$-norms of these discs are uniformly bounded  with respect to $t$, then then one can find $C> 0$ and $\tau> 0$ independent of  $t$ for the upper bound of  $\parallel\overline\partial_J (F \circ h_t) \parallel$.

\subsection{One-dimensional case} 

Recall some boundary properties of subsolutions of the $\overline\partial$-operator  in the unit disc.

Denote by $W^{k,p}(\D)$  the usual Sobolev classes of functions admitting generalized partial derivatives up to the order $k$ in $L^p(\D)$ (in fact we need only the case $k=0$ and $k=1$). In particular  $W^{0,p}(\D) = L^p(\D)$. We will always assume that $p > 2$.




Denote also by 
$\parallel f \parallel_\infty = \sup_\D \vert f \vert$
the usual $\sup$-norm on the space $L^\infty(\D)$ of complex functions bounded on $\D$.

\begin{lemma}
\label{SchwarzLemma}
Let $f \in L^\infty(\D)$ and $f_{\overline\zeta} \in L^p(\D)$ for some $p > 2$. Then 
\begin{itemize}
\item[(a)]  $f$ admits a non-tangential limit at almost every point $\zeta \in b\D$.
\item[(b)] if $f$ admits a limit along a  curve in $\D$ approaching $b\D$ non-tangentially at a boundary point $e^{i\theta} \in b\D$, then $f$ admits a non-tangential limit at $e^{i\theta}$.
\item[(c)] for each positive $r < 1$ there exists a constant $C = C(r) > 0$ (independent of $f$) such that for every $\zeta_j \in r\D$, $j=1,2$ one has 
\begin{eqnarray}
\label{SchwarzIn}
\vert f(\zeta_1)  - f(\zeta_2)\vert \le C (\parallel f \parallel_\infty + \parallel f_{\overline\zeta} \parallel_{L^p(\D)} ) \vert \zeta_1 - \zeta_2\vert ^{1-2/p}
\end{eqnarray}
\end{itemize}
\end{lemma}
The proof is contained in \cite{Su2}.

Sometimes it is convenient to apply the part (c) of Lemma on the disc $\rho \D$ with $\rho > 0$. Let $g \in L^\infty(\rho\D)$ and $g_{\overline\zeta} \in L^p(\rho\D)$. The function $f(\zeta):= g(\rho\zeta)$ satisfies the assumptions of Lemma \ref{SchwarzLemma} on $\D$. Let $0 < \alpha < \rho$ and let  $\vert \tau_j \vert < \alpha$, $j = 1,2$. Set $\zeta_j = \tau_j/\rho$.   Then $\vert \zeta_j \vert < r= \alpha/\rho < 1$, $j = 1,2$.
 Applying (c) Lemma \ref{SchwarzLemma} to $f$   we obtain:

\begin{eqnarray}
\label{Schwarzln1}
\vert g(\tau_1)  - g(\tau_2)\vert \le (C(r)/\rho^{1-2/p}) (\parallel g \parallel_\infty + \rho \parallel g_{\overline\zeta} \parallel_{L^p(\rho\D)} ) \vert \tau_1 - \tau_2\vert ^{1-2/p}
\end{eqnarray}
Note that the constant $C = C(r) = C(\alpha/\rho)$ depends only on the quotient $r = \alpha/\rho < 1$. If $r$ is separated from $1$, the value of  $C$ is fixed.

\section{Main results }

First we introduce an almost complex analog of an admissible approach which is classical in  the case of $\C^n$, see \cite{St,Ch}.

Let $\Omega$ be a smoothly bounded domain in an almost complex manifold $(M,J)$. Notice that any domain with boundary of class $C^2$ satisfies all assumptions imposed below. Fix a hermitian metric on $M$ compatible with $J$; a choice of such metric will not affect our results since it changes only constant factors in estimates. We measure all distances and norms with  respect to the choosen metric.

Let $p \in b\Omega$ be a boundary point. {\it A non-tangential approach} to $b\Omega$ at $p$ can be defined as the limit along the sets
\begin{eqnarray}
\label{NonTan1}
C_\alpha(p) = \{ q \in \Omega: dist(q,p) < \alpha \delta_p(q) \}, \,\,\, \alpha > 1.
\end{eqnarray}
Here $\delta_p(q)$ denotes the minimum of distances from $q$ to the tangent plane $T_p(b\Omega)$ 
and to $b\Omega$.

We need to define a wider class of regions. {\it An admissible approach}  to $b\Omega$ at $p$ is defined as the limit along the sets

\begin{eqnarray}
\label{Ad1}
A_{\alpha, \varepsilon}(p) = \{ q \in \Omega: d_p(q) < (1+\alpha)\delta_p(q), \, dist(p,q)^2 < \alpha \delta_p^{1+\varepsilon}(q) \}, \,\, \alpha > 0, \,\,\varepsilon > 0.
\end{eqnarray}
Here $d_p(q)$ denotes the distance from $q$ to {\it the holomorphic tangent space} $H_p(b\Omega) = T_p(b\Omega) \cap J T_p(b\Omega)$. Similarly to the classical case of $\C^n$, an admissible region approaches $b\Omega$ transversally in the normal direction and can be tangent in the directions of the holomorphic tangent space.

\begin{definition}
A function $F:\Omega \to \C$ has an admissible limit $L$ at $p \in b\Omega$ if 
$$\lim_{A_{\alpha,\varepsilon}(p) \ni q} F(q) = L \,\, \mbox{for all} \,\,  \alpha, \varepsilon > 0.$$
\end{definition}

Next we need the following notion.

\begin{definition}
Let $\Omega$ be a smoothly bounded domain in an almost complex manifold $(M,J)$ of complex dimension $n$ and $p \in b\Omega$ be a boundary point. A real curve 
$\gamma: [0,1[ \to \Omega$  of class $C^1([0,1])$  is called an admissible $p$-curve if $\gamma(1) = p$ and $\gamma$ is transverse to the tangent space $T_p(b\Omega)$ (i.e. the tangent vector of $\gamma$ at $p$ is not contained in $T_p(b\Omega)$).
\end{definition}

\begin{definition}
A function $F$ defined on $\Omega$ has a limit $L \in \C$ along an admissible $p$-curve $\gamma$ if there exists  $\lim_{t \to 1} (F \circ \gamma)(t) = L$.
\end{definition}


Our first main result is the following analog of the Chirka - Lindel\"of principle \cite{Ch}.

\begin{thm}
\label{Thm1}
Let $\Omega$ be a smoothly bounded  domain in an almost complex manifold $(M,J)$ of complex dimension $n$. Suppose that a complex function  $F \in L^{\infty}(\Omega)$   is a  
$\overline\partial_J$ -subsolution (in the sense of Definition \ref{DefSub}) on $\Omega$. If $F$ has a limit  along an admissible $p$-curve for some $p \in b\Omega$, then $F$ has an admissible limit  at $p$.
\end{thm}
A similar result is obtained in \cite{Su2} under considerably stronger assumptions. First, the domain $\Omega$ in \cite{Su2} is supposed strictly pseudoconvex. Second, it is assumed that an admissible curve $\gamma$ is contained in some pseudoholomorphic disc. Finally, the assumption (\ref{AlHol}) is replaced there by the stronger condition of  boundedness of 
$\overline \partial_J  F(z)$ on $\Omega$.



\bigskip

As an application of the Chirka -Lindel\"of principle we establish a Fatou type results for $\overline\partial_J$-subsolutions. For holomorphic functions in $\C^n$ the first versions of the Fatou theorem are due to E.Stein \cite{St},  E.Chirka \cite{Ch}, F. Forstneri\v c \cite{Fo} and  A.Sadullaev \cite{Sa1}. Our approach is inspired by \cite{Su2}.

We will deal with  some standard  classes of real submanifolds of an almost complex manifold. A submanifold $E$ of an almost complex $n$-dimensional $(M,J)$ is called 
{\it totally real} if at every point $p \in E$ the tangent space $T_pE$  does not contain non-trivial complex vectors that is $T_pE \cap JT_pE = \{ 0 \}$. This is well-known that the (real)  dimension of a totally real submanifold of $M$ is not bigger than $n$; we will consider in this paper only $n$-dimensional totally real submanifolds that is the case of maximal dimension. A real submanifold $N$ of $(M,J)$ is called {\it generic} if the complex span of $T_pN$ is equal to the whole  $T_pM$ for each point $p \in N$. A real $n$-dimensional submanifold of $(M,J)$ is generic if and only if it is totally real.  





\bigskip

Our second main result here is the following

\begin{thm}
\label{Thm2}
Let $E$ be a generic submanifold of the boundary $b\Omega$ of a smoothly bounded domain $\Omega$ in an almost complex manifold $(M,J)$ of complex dimension $n$. Suppose that a complex function $F \in L^{\infty}(\Omega)$ is a $\overline\partial_J F(p)$-subsolution   on $\Omega$. Then $F$ has an admissible limit  at almost every point of $E$.
\end{thm}
Note that the Hausdorff $n$-meausure on $E$ here is defined with respect to any metric on $M$; the condition to be a subset of measuro zero in $E$ is independent of such a choice.
A similar result also is obtained in \cite{Su2} under considerably stronger assumptions discussed above: the domain $\Omega$ in \cite{Su2} is supposed strictly pseudoconvex and  the assumption (\ref{AlHol}) is replaced there by the stronger condition of  boundedness of 
$\overline \partial_J  F(z)$ on $\Omega$.

\section{Proofs}
Assume that we are in the setting of Theorem \ref{Thm1}.  First we need the following

\begin{lemma}
\label{LindTang}
Let $F$ satisfies assumptions of Theorem \ref{Thm1}. If $F$ has a limit  along a $p$-admissible curve $\gamma_1$ at $p \in E$, then $F$ has the same limit along each admissible curve in  $\Omega$ tangent to $\gamma_1$ at $p$.
\end{lemma}
\proof  Let $\gamma_2$ be another $p$-admissible curve and such that $\gamma_1$ and $\gamma_2$ have the same tangent line at $p$. Without loss of generality asssume that $p = 0$ (in local coordinates). Denote by $\rho$ a local defining function of $\Omega$. 


It follows by the Nijenhuis-Woolf theorem that there exists a family $z_t( \zeta): \D \to \C^n$,  of embedded $J$-holomorphic discs  near the origin in $\C^n$ satisfying the following properties:
\begin{itemize}
\item[(i)] the family $z_t$ is smooth on $\overline \D \times [0,1]$ 
\item[(ii)] for every $t \in [0,1]$  the disc $z_t$ transversally intersects each curve $\gamma_j$  at a unique point coresponding to some  parameter value $\zeta_j(t) \in \D$ ,  $j = 1,2$. In other words $\gamma_j(t) = z_{t}( \zeta_j(t))$. Furthermore, $\zeta_1(t) = 0$, i.e. this point is the center of the disc $z_t$.
\end{itemize}
In the case of the standard complex structure each such a disc is simply  an open piece (suitably parametrized) of a complex line intersecting transversally  the both of curves $\gamma_j$. Recall that the curves are embedded near the origin and tangent at the origin so such a family of complex lines obviously exist. The $J$-holomorphic discs are obtained  from this family of lines by a small deformation described in the proof of the Nijenhuis-Woolf theorem in Section 2. Note that for $t=1$ the disc $z_1$ intersects transversally the both curves $\gamma_j$ at the same point $\gamma_j(1) = p$.

Furthermore, because of the condition (i),  the  compositions $F \circ z_t$ have $\overline\zeta$- derivatives  of class $L^p$ on their domains of definitions, for each $p > 2$ close enough to $2$. Moreover,  their $L^p$ norms are  bounded  on $\D$ uniformly with respect to  $t$. Indeed, it follows by the Chain Rule and (\ref{holomorphy}) that
$$(F \circ z)_{\overline\zeta} = (F_{\overline z} + F_zA){\overline z}_{\overline\zeta}$$
and now we use the assumption that  $\overline\partial_J F(z)$ has the growth of order $dist(z,b\Omega)^{-1/2+ \tau}$, $\tau > 0$.


Since the curves $\gamma_j$ are tangent at the origin, we have

\begin{eqnarray}
\label{radius}
\vert \zeta_2(t) \vert = o(1-t)
\end{eqnarray}
as $t  \to 1$.

The curve $\gamma_1$ is admissible, so we have
$$dist(\gamma_1(t), b\Omega) = O(1-t)$$
as $t \to 1$.   Hence, there exists $\rho(t) = O(1-t)$ as $t \to 1$ such that $z_t(\rho(t)\D)$ is contained in $\Omega$.  Applying (\ref{Schwarzln1}) to the composition $f:= F \circ z_t(\zeta)$  on the disc $\rho(t)\D$, we obtain
(fixing $r > 0$)

\begin{eqnarray}
\label{Schwarzln2}
\vert f(0)  - f(\zeta_2(t))\vert \le (C/O(1-t)^{1-2/p}) (\parallel f \parallel_\infty + O(1-t) \parallel f_{\overline\zeta} \parallel_{L^p} ) o((1-t)^{1-2/p} \to 0
\end{eqnarray}
as $t  \to 1$. Note that by (\ref{radius}) for every $t$ the point $\zeta_2(t)$ is contained in $(1/2)\rho(t)\D$; hence, the constant $C$ is independent of $t$ (see remark after (\ref{Schwarzln1}) ).
This concludes the proof of Lemma.

\bigskip

Now we continue the proof of Theorem. Consider first the special case where our almost complex manifold $M$ coincides with $\C^n$ and 
the almost complex structure $J$ coincides with $J_{st}$.  

It suffices to consider the  case where $p = 0$. Furthermore, after a linear change of coordinates we have the defining function $\rho$ of $\Omega$ has the form

\begin{eqnarray}
\label{DefFun}
\rho(z) = y_n + o(\vert z \vert)
\end{eqnarray}

In particular, the holomorphic tangent space $H_0(b\Omega)$ has the form
\begin{eqnarray}
\label{HolSp1}
H_0(b\Omega) = \{ z_n = 0 \}
\end{eqnarray}
Without loss of generality we employ the usual Euclidean distance.


We have  $T_0(b\Omega_0) = \{ y_n = 0 \}$.  Note that
$$dist(z,b\Omega_0) \sim \vert \rho(z) \vert \le \vert y_2 \vert = dist(z,T_0(b\Omega_0)).$$ 
 Hence we can assume $\delta_0(z) = \vert\rho(z)\vert$. 
Since $dist(z,H_0(b\Omega_0)) = \vert z_n \vert$, for each $\alpha > 0$ the admissible regions $A_{\alpha,\varepsilon}(0)$  from (\ref{Ad1}) are defined by the conditions

\begin{eqnarray}
\label{angle3}
\vert z_n \vert < (1 + \alpha) \vert \rho(z) \vert
\end{eqnarray}
and 
\begin{eqnarray}
\label{angle4}
 \vert z \vert^2 < \alpha\vert \rho(z) \vert^{1+\varepsilon}
\end{eqnarray}

After an additional linear change of coordinates (which preserves the previous setting) one can assume that the tangent line $T_0(\gamma)$ is contained in the coordinate complex line 
$L_n = (0,...,0,z_n)$, $z_n \in \C$. By Lemma \ref{LindTang} the function $F$ admit the limit along the ray in $L_n$ which is tangent to $T_0(\gamma)$.





The intersection of the complex normal plane $L_n$ with  $\Omega$  is the plane domain $\Pi:= \{ z: z_1=...=z_{n-1} = 0,  y_n + o(y_n) < 0 \}$ and the first inequality (\ref{angle3}) defines a non-tangential region there (which tends to this half-plane when $\alpha$ increases). 

Fix a point $(0,...,0,z_n^0)$ which satisfies (\ref{angle3}). Fix a unit vector $v \in H_0(b\Omega)$ of the form  $v = (v_1,...,v_{n-1}, 0)$.  Consider a complex line through the point $(0,...,0,z^0_n)$ in the direction $v$:

\begin{eqnarray}
\label{Heisn2}
f(v,z^0_n): \C \ni \zeta \mapsto (\zeta v, z^0_n)
\end{eqnarray}
which is  parallel to $H_0(b\Omega_0)$. A simple calculation shows that the second assumption  (\ref{angle4}) is equivalent to the fact that $f(v,z_n^0)(r\D) \subset A_{\alpha,\varepsilon}(0)$ when
\begin{eqnarray}
\label{Heisn3}
r  \sim \vert y_n^0\vert^{1/2+ \varepsilon}
\end{eqnarray}
Clearly, this family of complex discs fill  the region $A_{\alpha,\varepsilon}(0)$ when $(0,...,0, z_n^0)$ satisfies the first condition (\ref{angle3}). Furtermore, since $\rho^0:= \vert \rho(0,...,0,z^0_n) \vert \sim \vert y_n^0 \vert$, the disc $f(v,z_n^0)(\rho^0\D)$ is contained in $\Omega$.


The restriction $F \circ L_n$ is a bounded function on $\Pi$ and $(F \circ L_n)_{\overline\zeta}$ is is of class $L^p$ with $p > 2$ close enogh to $2$. Furthermore, $F \circ L_n$ admits a limit $L$ along some ray in $\Pi$ with vettex at $0$. By (b) Lemma \ref{SchwarzLemma} the function $F \circ L_n$ admits the limit $L$ along any non-tangential region in $\Pi$.
Let now $z \in A_{\alpha,\varepsilon}(0)$. Then there exists a unit vector $v \in H_0(b\Omega)$ and a point $z_n^0$ in the non-tangential region on $\Pi$ such that the disc $f(v,z^0_n)$  contains the point $z$ that is $z = f(v,z_n^0)(\zeta)$ for some $\zeta$ with $\vert \zeta \vert \le C \vert y_n^0 \vert^{1/2}$. Since also $f(v,z_n^0)(0) = z_n^0$, 
by (\ref{Schwarzln1}) we have the estimate:
$$\vert F(z) - F(0,...,0,,z_n^0) \vert = \vert (F \circ f(v,z_n^0))(\zeta) - (F \circ f(v,z_n^0))(0) \vert  \le C  \vert y_n^0 \vert^\tau$$
with $\tau = \varepsilon ( 1 - 2/p) > 0$. Note that we apply (\ref{Schwarzln1}) on a disc $\rho^0 \D$ and use (\ref{Heisn3}) because $z \in A_{\alpha,\varepsilon}(0)$.
Since $F(0,...,0,z_n) \to L$ as $y_n^0 \to 0$, we conclude that $F(z) \to L$.

The case of a general almost complex structure $J$ follows by the same argument using a slight deformation transforming the above $J_{st}$-holomorphic discs to $J$-holomorphic discs.  Such a deformation is always possible by the Nijenhuis-Woolf theorem and is continuous in any $C^k$ norm. Thus, it changes only constants in estimates and the above argument literally goes through. This proves Theorem \ref{Thm1}.

\bigskip

Now the proof of Theorem \ref{Thm2} follows exactly as in \cite{Su2} using Theorem \ref{Thm1}.

\end{document}